\documentclass{amsart}

\usepackage{fancyhdr}
\usepackage{lastpage}
\usepackage{yhmath}

\newcommand{\bdis}{\begin{displaymath}}
\newcommand{\edis}{\end{displaymath}}
\newcommand{\be}{\begin{equation}}
\newcommand{\ee}{\end{equation}}

\newcommand{\mbb}{\mathbb}
\newcommand{\mcal}{\mathcal}

\newcommand{\vp}{\varphi}
\newcommand{\vt}{\vartheta}

\newcommand{\zf}{\zeta\left(\frac{1}{2}+it\right)}

\DeclareMathOperator{\re}{Re}

\DeclareMathOperator{\Var}{Var}

\theoremstyle{definition}

\theoremstyle{remark}
\newtheorem{remark}[]{Remark}

\newtheorem*{mydef1}{{\bf Theorem}}

\newtheorem*{mydef2}{{\bf Definition}}

\newtheorem*{mydef4}{{\bf Corollary}}

\newtheorem*{mydef51}{{\bf Lemma 1}}

\newtheorem*{mydef52}{{\bf Lemma 2}}

\newtheorem*{mydef53}{{\bf Lemma 3}}

\numberwithin{equation}{section}

\begin{document}

\title{On statistical arc length of the Riemann $Z(t)$-curve}

\author{Jan Moser}

\address{Department of Mathematical Analysis and Numerical Mathematics, Comenius University, Mlynska Dolina M105, 842 48 Bratislava, SLOVAKIA}

\email{jan.mozer@fmph.uniba.sk}


\begin{abstract}
In this paper we study certain stochastic process that is generated by the Riemann-Siegel formula. Further, we construct corresponding statistical
model by a way similar to those used in telecommunication. We define statistical arc length of the Riemann $Z(t)$-curve in this model and obtain
an asymptotic formula for that length. This paper is English remake of our work of reference \cite{4}.
\end{abstract}

\maketitle

\section{Introduction}

\subsection{}

In the paper \cite{5} we have studied the following integral
\be \label{1.1}
\int_T^{T+H}\sqrt{1+\{ Z'(t)\}^2}{\rm d}t=L(T,H)=L,
\ee
i.e. the arc length of the Riemann curve
\bdis
y=Z(t),\ t\in [T,T+H],\ T\to\infty,
\edis
where (see \cite{7}, pp. 79, 329)
\bdis
\begin{split}
& Z(t)=e^{i\vt(t)}\zf , \\
& \vt(t)=-\frac t2\ln\pi +\mbox{Im}\ln\Gamma\left(\frac 14+i\frac t2\right).
\end{split}
\edis
Next, we will denote the roots of the equations
\bdis
Z(t)=0,\quad Z'(t)=0
\edis
by the symbols
\bdis
\{\gamma\},\quad \{ t_0\};\quad t_0\not=\gamma
\edis
correspondingly.

\begin{remark}
On the Riemann hypothesis, the points of the sequences $\{\gamma\}$ and $\{ t_0\}$ are separated each from other (see
\cite{2}, Corollary 3), i.e. in this case we have
\bdis
\gamma'<t_0<\gamma'' ,
\edis
where $\gamma',\gamma''$ are neighboring points of the sequence $\{\gamma\}$. Of course, $Z(t_0)$ is an extremal value
of the function $Z(t)$ in some neighborhood of the point $t_0$.
\end{remark}

Namely, we have proved (see \cite{5}) for integral (\ref{1.1}) the following formula
\be \label{1.2}
\begin{split}
& \int_T^{T+H}\sqrt{1+\{ Z'(t)\}^2}{\rm d}t=2\sum_{T\leq t_0\leq T+H}|Z(t_0)|+
 \Theta H+\mcal{O}\left(T^{\frac{A}{\ln\ln T}}\right), \\
 &  T\to\infty,\ \Theta=\Theta(T,H)\in (0,1),\ H=T^{\epsilon}
\end{split}
\ee
for every fixed $\epsilon>0$.

The proof of the formula (\ref{1.2}) was based on the following
\be \label{1.3}
\begin{split}
 & Z'(t)=-2\sum_{n<P}\frac{1}{\sqrt{n}}\ln\frac{P}{n}\sin\{\vt(t)-t\ln n\}+ \\
 & + \mcal{O}(T^{-1/4}\ln T),\ t\in [T,T+U],\ U\in (0,\sqrt{T}],
\end{split}
\ee
that was obtained in our paper \cite{3} as a variant of the Riemann-Siegel formula. One sight on formulae
(\ref{1.1})--(\ref{1.3}) is sufficient for insight of the difficulty of the problem about the asymptotic
formula of the arc length of the Riemann curve.

\subsection{}

In connection with formulae (\ref{1.1})--(\ref{1.3}) we give the following

\begin{remark}
The main difficulty lies probably in that we are absent of something like famous Dirac procedure
\bdis
\sqrt{m^2c^2+p_1^2+p_2^2+p_3^2}=\alpha_1p_1+\alpha_2p_2+\alpha_3p_3+\beta,
\edis
i.e. extracting of the square root of the operator - see \cite{1}, p. 255, equations (5), (7). This was the
way towards Dirac's relativistic wave equation for the electron - the fundamental law of the quantum
mechanics.
\end{remark}

\subsection{}

Next, let us remind, in connection with (\ref{1.1})--(\ref{1.3}), that in our recent paper we have proved
the following result on metamorphoses: there is an infinite set of elements
\be \label{1.4}
\{ \alpha_0(T),\alpha_1(T),\dots,\alpha_k(T)\},\ T\in (T_0,+\infty),\ T_0>0
\ee
where $T_0$ is sufficiently big, such that
\be \label{1.5}
\begin{split}
 & \prod_{r=1}^k
 \left|
 \sum_{n\leq\tau(\alpha_r)}\frac{2}{\sqrt{n}}
 \cos\{\vt(\alpha_r)-\alpha_r\ln n\}+\mcal{O}(\alpha_r^{-1/4})
 \right| \sim \\
 & \sim
 \sqrt{\frac{\Lambda}
 {\left|\sum_{n\leq\tau(\alpha_0)}\frac{2}{\sqrt{n}}
 \cos\{\vt(\alpha_0)-\alpha_0\ln n\}+\mcal{O}(\alpha_0^{-1/4})\right|}} , \\
 & T\to\infty ,
\end{split}
\ee
where
\bdis
\Lambda=\sqrt{2\pi}\frac{H}{H_k}\ln^kT,\ \tau(t)=\sqrt{\frac{t}{2\pi}},\ k=1,\dots,k_0,\ k_0\in\mbb{N},
\edis
i.e. to the infinite subset
\bdis
\{\alpha_1(T),\dots,\alpha_k(T)\}
\edis
an infinite set of metamorphoses of the multiform on the left-hand side of (\ref{1.5}) into quite
distinct form on the right-hand side of (\ref{1.5}) corresponds.

Now we rewrite the formula (\ref{1.5}) in the following form
\be \label{1.6}
\begin{split}
 & \sqrt{\frac{1}{\Lambda}
 \left|\sum_{n\leq\tau(\alpha_0)}\frac{2}{\sqrt{n}}
 \cos\{\vt(\alpha_0)-\alpha_0\ln n\}+\mcal{O}(\alpha_0^{-1/4})\right|}\sim \\
 & \sim
 \prod_{r=1}^k
 \left|
 \sum_{n\leq\tau(\alpha_r)}\frac{2}{\sqrt{n}}
 \cos\{\vt(\alpha_r)-\alpha_r\ln n\}+\mcal{O}(\alpha_r^{-1/4})
 \right|^{-1}.
\end{split}
\ee

\begin{remark}
We see, in connection with (\ref{1.1})--(\ref{1.3}) and the Remark 2, that the square root of the
weighted nonlinear monoform on the left-hand side of (\ref{1.6}) is asymptotically expressed as the rational
function of the multiform for infinite set of elements (\ref{1.4}).
\end{remark}

\subsection{}

In this paper, we introduce new method to study the main integral (\ref{1.1}). Namely, we use certain stochastic
process generated by the formula (\ref{1.3}). We use simple properties of this process together with the
Ljapunov's central limit theorem for a construction of an integral that we call
\emph{statistical arc length of the Riemann $Z(t)$-curve}. Consequently, an asymptotic formula for this
integral is obtained.

\section{Tranformation of the formula (\ref{1.1})}

In the formula (see (\ref{1.1}))
\bdis
L=L(T,U)=\int_T^{T+U}\sqrt{1+\{ Z'(t)\}^2}{\rm d}t
\edis
we put
\bdis
Z'(t)=Z_1(t)+R_1(t),\ t\in [T,T+U],
\edis
where (see (\ref{1.3}))
\bdis
\begin{split}
& Z_1(t)=2\sum_{n<P}\frac{1}{\sqrt{n}}\ln\frac Pn\cos\{\vt(t)-t\ln n+\pi/2\}, \\
& R_1(t)=\mcal{O}(T^{-1/4}\ln T).
\end{split}
\edis
Since
\bdis
\ln\frac Pn=\mcal{O}(\ln P),\ \ln\frac{P}{n-1}>\ln\frac Pn,\ 2\leq n<P,
\edis
then, similarly to \cite{7}, pp. 92, 93, we obtain the following estimate
\bdis
Z_1(t)=\mcal{O}(T^{1/6}\ln^2T),\ t\in [T,T+U].
\edis
Hence,
\bdis
\{Z'(t)\}^2=\{Z_1(t)\}^2+\mcal{O}(T^{-1/12}\ln^3T),
\edis
and
\bdis
\begin{split}
& \sqrt{1+\{Z'(t)\}^2}=\sqrt{1+\{ Z_1(t)\}^2}\left\{ 1+\mcal{O}\left(\frac{T^{-1/12}\ln^3T}{1+\{ Z_1(t)\}^2}\right)\right\}= \\
& = \sqrt{1+\{ Z_1(t)\}^2}+\mcal{O}(T^{-1/12}\ln^3T).
\end{split}
\edis
Consequently, we have the following

\begin{mydef51}
\be \label{2.1}
L=L(T,U)=\int_T^{T+U}\sqrt{1+\{ Z_1(t)\}^2}{\rm d}t+\mcal{O}(T^{-1/12}\ln^3T),
\ee
\be \label{2.2}
Z_1(t)=2\sum_{n<P}\frac{1}{\sqrt{n}}\ln\frac Pn\cos\{\vt(t)-t\ln n+\pi/2\}, \quad U\in (0,\sqrt{T}].
\ee
\end{mydef51}

\section{Definition of the statistical arc length of the Riemann curve and Theorem}

\subsection{}

Let
\bdis
\vp_n\in [-\pi,\pi],\ n<P
\edis
be the system of independent random variables each of them uniformly distributed on the segment $[-\pi,\pi]$.
Putting these $\vp_n$ into the arguments of cosine-functions in (\ref{2.2}) one obtains the following
stochastic process
\be \label{3.1}
\begin{split}
& \Phi_1(t)=\Phi_1(t,\vp_1,\vp_2,\dots)= \\
& = 2\sum_{n<P}\frac{1}{\sqrt{n}}\ln\frac Pn\cos\{ a(t,n)+\vp_n\},
\end{split}
\ee
where
\bdis
\begin{split}
& a(t,n)=\vt(t)-t\ln n+\frac \pi2 , \\
& t\in [T,T+U],\ U\in (0,\sqrt{T}].
\end{split}
\edis

\begin{remark}
The following is true: every realization
\bdis
\Phi_1(t,\bar{\vp}_1,\bar{\vp}_2,\dots),\ t\in [T,T+U]
\edis
of the stochastic process (\ref{3.1}) is a continuous function of the variable $t$ for every admissible vector
\bdis
\bar{\vp}=(\bar{\vp}_1,\bar{\vp}_2,\dots)
\edis
and. consequently, there is the Riemann integral
\bdis
\int_T^{T+U}\Phi_1(t,\bar{\vp}){\rm d}t,\quad \forall \bar{\vp}.
\edis
\end{remark}

We also use the one-dimensional Gaussian density of probability
\be \label{3.2}
w_\infty(\Phi_1)=\frac{1}{2}\frac{1}{\ln^{3/2}P}e^{-\beta (\Phi_1)^2},\quad \beta=\frac{3}{4\ln^3P}.
\ee

\subsection{}

Now, we define the following stochastic process
\be \label{3.3}
\Phi_2=\Phi_2(T,U,\vp)=\int_T^{T+U}\sqrt{1+\{\Phi_1(t)\}^2}{\rm d}t,\ U\in (0,\sqrt{T}],
\ee
for every sufficiently big $T$. Expectation value of this process is given by the formula
\be \label{3.4}
E(\Phi_2)=\int_{T}^{T+U} E(\sqrt{1+(\Phi_1)^2}){\rm d}t.
\ee
Next, we define the following asymptotic expectations
\be \label{3.5}
E_\infty(\sqrt{1+(\Phi_1)^2})=\int_{-\infty}^\infty \sqrt{1+(\Phi_1)^2}w_\infty(\Phi_1){\rm d}\Phi_1,
\ee
and, consequently, (see (\ref{3.3})--(\ref{3.5}))
\be \label{3.6}
E_\infty(\Phi_2)=\int_T^{T+U} E_\infty(\sqrt{1+(\Phi_1)^2}){\rm d}t.
\ee
Thus, the comparison of the formulae (\ref{2.2}), (\ref{3.1}) and (\ref{2.1}), (\ref{3.6}) leads us
to the following

\begin{mydef2}
\be \label{3.7}
\langle L(t,U) \rangle|_S=E_\infty(\Phi_2),
\ee
where
\bdis
\langle L(t,U) \rangle|_S
\edis
is the statistical arc length of the Riemann $Z(t)$-curve.
\end{mydef2}

Consequently, we notice that the basis of constructed statistical model lies in the following (purely) mathematical

\begin{mydef1}
\be \label{3.8}
E_\infty(\sqrt{1+(\Phi_1)^2})\sim \frac{1}{\sqrt{6\pi}}\ln^{3/2}T,\ T\to\infty.
\ee
\end{mydef1}

Since (see (\ref{3.6}), (\ref{3.8}))
\be \label{3.9}
E_\infty(\Phi_2)\sim\frac{1}{\sqrt{6\pi}} U\ln^{3/2}T,\ T\to\infty ,
\ee
then we obtain from (\ref{3.9}) the following

\begin{mydef4}
\bdis
\langle L(t,U) \rangle|_S\sim\frac{1}{\sqrt{6\pi}} U\ln^{3/2}T,\ U\in (0,\sqrt{T}],\ T\to\infty.
\edis
\end{mydef4}

\section{Statistical considerations about the Gaussian distribution}

\subsection{}

The following lemma holds true

\begin{mydef52}
\be \label{4.1}
E(\Phi_1)=0,\ \Var(\Phi_1)\sim \frac 23\ln^3P,\ T\to\infty .
\ee
\end{mydef52}

\begin{proof}
We have (see (\ref{3.1}))
\bdis
\Phi_1(t;\vp)=\sum_{n<P} X_n,\ X_n=\frac{2}{\sqrt{n}}\ln\frac Pn\cos\{ \vp_n+a\}.
\edis
Since $\vp_n$ are independent and uniformly distributed then we obtain that
\bdis
E(X_n)=0,\ \Var(X_n)=\frac 2n\ln^2\frac Pn
\edis
and
\bdis
E(\Phi_1)=0,\ \Var(\Phi_1)=2\sum_{n<P} \frac 1n\ln^2\frac Pn.
\edis
Next, for $\Var$ we use the Euler-MacLaurin summation formula (comp. \cite{7}, p. 13) for the function
\bdis
f(x)=\frac 1x\ln^2\frac Px,\ x\in [1,P].
\edis
Since
\bdis
f'(x)=\mcal{O}\left(\frac{\ln^2P}{x^2}\right),\ f(1)=\mcal{O}(\ln^2P),\ f(P)=0,
\edis
then
\bdis
\begin{split}
 & \sum_{n<P}f(n)=\int_1^P\ln^2\frac Px \frac{{\rm d}x}{x}+
 \mcal{O}\left(\ln^2P\int_1^P\frac{{\rm d}x}{x^2}\right)+\mcal{O}(\ln^2P)= \\
 & = \frac 13\ln^3P+\mcal{O}(\ln^2P),
\end{split}
\edis
and, consequently,
\bdis
\Var(\Phi_1)\sim \frac 23\ln^3P,\ T\to\infty.
\edis
\end{proof}

\subsection{}

Next, we give the following

\begin{mydef53}
It is true that in our case the Ljapunov's condition in the central limit theorem is fulfilled.
\end{mydef53}

\begin{proof}
Since
\bdis
\begin{split}
 & E(|X_n-E(X_n)|^3)=E(|X_n|^3)= \frac{8}{n^{3/2}}\ln^3P\frac{1}{2\pi}
 \int_{-\pi}^\pi |\cos^3(\vp_n+a)|{\rm d}\vp_n< \\
 & < \frac{8}{n^{3/2}}\ln^3\frac Pn,
\end{split}
\edis
then we have the following estimate for the third absolute central moment
\bdis
\begin{split}
 & \sum_{n<P}E(|X_n|^3)\leq 8\sum_{n<P}\frac{1}{n^{3/2}}\ln^3\frac Pn<8\ln^3P
 \sum_{n=1}^\infty \frac{1}{n^{3/2}}=A\ln^3P.
\end{split}
\edis
Now, (see (\ref{4.1}))
\bdis
B_P=\sqrt{\Var(\Phi_1)}\sim \sqrt{\frac 23}\ln^{3/2}P
\edis
and, consequently,
\be \label{4.2}
\frac{1}{B_P^{3}}\sum_{n<P}E(|X_n|^3)<2A\frac{\ln^3P}{\ln^{9/2}P}\xrightarrow{T\to\infty}0,
\ee
i.e. the Ljapunov's condition is fulfilled.
\end{proof}

\begin{remark}
Our choice of the Gaussian asymptotics in (\ref{3.2}) is based upon the results (\ref{4.1}) and (\ref{4.2}),
respectively.
\end{remark}

\section{Proof of Theorem}

We put (see (\ref{3.2}) -- (\ref{3.5}))
\be \label{5.1}
E_\infty(\sqrt{1+(\Phi_1)^2})=\sqrt{\frac{3}{\pi}}\frac{1}{\ln^{3/2}P}F(\beta),
\ee
where
\be \label{5.2}
F(\beta)=\int_0^\infty \sqrt{1+x^2}e^{-\beta x^2}{\rm d}x,\ \beta=\frac{3}{4\ln^3P}.
\ee
We need to obtain an asymptotic formula for $F(\beta)$ for small values of $\beta$ (i.e. for
large values of $T$).

First of all, we have that
\bdis
\begin{split}
 & F(\beta)=\int_0^\infty \cosh^2te^{-\beta\sinh^2t}{\rm d}t=
 \frac 12 e^{\beta/2}\int_0^\infty (\cosh 2t+1)e^{-\frac \beta2 \cosh 2t}{\rm d}t= \\
 & = \frac 14 e^{\beta/2}\int_0^\infty (\cosh t+1)e^{-\frac \beta2 \cosh t}{\rm d}t.
\end{split}
\edis
Now, we use the Schl\" afli's integral (see \cite{8}, p. 81)
\bdis
\int_0^\infty e^{-z\cosh t}\cosh\nu t{\rm d}t=K_\nu(t),\ \re\{z\}>0,
\edis
for the modified Bessel's function of the second kind, and then we obtain
\be \label{5.3}
F(\beta)=\frac 14 e^{\beta/2}\left\{ K_0\left(\frac \beta2\right)+K_1\left(\frac \beta2\right)\right\}.
\ee
Further, we use the following representations of the $K_\nu$-functions (see \cite{8}, p. 80)
\bdis
\begin{split}
 & K_0(z)=I_0(z)\ln\frac z2+\sum_{m=1}^\infty \frac{z^{2m}}{2^m(m!)^2}\psi(m+1), \\
 & K_n(z)=\frac 12\sum_{m=0}^{n-1} (-1)^m\frac{(n-m-1)!}{m!\left(\frac z2\right)^{n-2m}}+ \\
 & + (-1)^{n+1}\sum_{m=0}^\infty \frac{\left(\frac z2\right)^{n+2m}}{m!(n+m)!}
 \left\{ \ln\frac z2-\frac 12\psi(m+1)-\frac 12\psi(n+m+1)\right\},
\end{split}
\edis
where (see \cite{8}, p. 77)
\bdis
I_0(z)=\sum_{m=0}^\infty \frac{1}{(m!)^2}\left(\frac z2\right)^{2m},
\edis
and (see \cite{9}, p. 241)
\bdis
\psi(x)=\frac{{\rm d}}{{\rm d}x}\ln\Gamma(x)=-c-\sum_{k=0}^\infty
\left(\frac{1}{k+x}-\frac{1}{k+1}\right) ,
\edis
where $c$ is the Euler's constant.

Since
\bdis
\psi(m+1),\psi(n+m+1)=\mcal{O}(1),
\edis
and (see (\ref{5.2}))
\bdis
\begin{split}
 & \frac \beta2=\frac{3}{8\ln^3P},\ \ln\frac\beta4=-3\ln\ln P+\mcal{O}(1), \\
 & I_0\left(\frac\beta2\right)=1+\mcal{O}\left(\frac{1}{\ln^6P}\right),
\end{split}
\edis
then
\bdis
\begin{split}
 & K_0\left(\frac\beta2\right)=3\ln\ln P+\mcal{O}(1), \\
 & K_1\left(\frac\beta2\right)=\frac 2\beta+\mcal{O}(\beta|\ln\beta|)=
 \frac 83\ln^3P+\mcal{O}\left(\frac{\ln\ln P}{\ln^3P}\right).
\end{split}
\edis
Consequently (see (\ref{5.3}))
\bdis
F(\beta)=\frac 23\ln^3P+\mcal{O}(\ln\ln P),
\edis
and from this (see (\ref{5.1})) the formula (\ref{3.8}) follows, of course,
\bdis
\ln^{3/2}P\sim\frac{1}{2\sqrt{2}}\ln^{3/2}T,\ T\to\infty.
\edis

\thanks{I would like to thank Michal Demetrian for helping me with the electronic version of this work.}

\end{document}